# Solid Angle of Conical Surfaces, Polyhedral Cones, and Intersecting Spherical Caps

*Oleg Mazonka, 2012*

**Abstract** *This paper presents formulae for calculation the solid angle of intersecting spherical caps, conical surfaces and polyhedral cones.*

## 1. Introduction

The problem of calculating solid angles appears in many areas of science and applied mathematics. While calculation of the solid angle of a right circular cone is a simple exercise, calculation of solid angles of arbitrary conical shapes is not trivial. In this paper I would like to show a few solutions for solid angle of conical shapes and as a special case the intersection of two cones.

## 2. Conical surfaces

### 2.1 Solid angle of conical surface

Let a conical surface be defined by a closed parametric curve in space with apex positioned at the origin of coordinate system. This parametric curve can be projected onto the unit sphere by renormalising the distance to each point on the curve. Let us parameterise a projected curve by the angle (being a linear distance on the surface of the unit sphere) $l$, so that the curve on the sphere is defined as $\vec{s}(l)$, and denote the first and the second derivatives with symbols $u$ and $\tau$:

$$\vec{\tau} = \frac{d}{dl}\vec{s} \qquad \vec{u} = \frac{d^2}{dl^2}\vec{s} \qquad\qquad 1$$

Since the curve is quite arbitrary it can be formed by a number of vertices (corners) connected by smooth (up to second derivative) arcs. Let $\delta_i$ be a turn angle of a corner $i$, so that $\delta_i$ and the inward angle are supplementary. $u$ and $\tau$ are not defined at corners. For each corner $i$ it is possible to define the tangent $\tau_i^-$ of one side of the corner and the tangent $\tau_i^+$ of the other side. Having both tangents, the angle $\delta_i$ can be calculated via trigonometric functions, such as:

$$\tan\delta_i = \frac{\sin\delta_i}{\cos\delta_i} = \frac{|\tau_i^- \times \tau_i^+|}{\tau_i^- \cdot \tau_i^+} \qquad 2$$

As will be proven later in this section, the solid angle of a conical surface can be calculated as

$$\boxed{\Omega = 2\pi - \sum_i \delta_i - \oint dl \sqrt{\vec{u}^2 - (\vec{s}\cdot\vec{u})^2}} \qquad 3$$

where the sum is taken over all corners $\delta_i$ and the line integral is taken along the closed curve excluding corners (in which $u$ is undefined).

A special case of the closed curve $\vec{s}$ consisting of segments of great circles gives the known formula (Ref [1]) for a spherical polygon:

$$\Omega = 2\pi - \sum_{i=1}^n \delta_i = \sum_{i=1}^n (\pi - \delta_i) - (n-2)\pi \qquad 4$$

because the integral along segments of great circles is equal to zero. In the absence of corners on the closed curve the second derivative $u$ is well-defined everywhere and hence the above equation can be written as:

$$\Omega = 2\pi - \oint dl \sqrt{\vec{u}^2 - (\vec{s}\cdot\vec{u})^2} \qquad 5$$

## 2.2  Right circular cone example

For example, a circle on the unit sphere with an apex angle $\theta$ can be parameterised as:

$$\vec{s}(\varphi) = (\sin\theta\cos\varphi, \sin\theta\sin\varphi, \cos\theta) \qquad 6$$

Then taking into account that $dl = d\varphi\sin\theta$

$$\vec{u} = \frac{d^2\vec{s}}{dl^2} = \frac{d^2\vec{s}}{d\varphi^2}\left(\frac{d\varphi}{dl}\right)^2 = \frac{1}{\sin^2\theta}(-\sin\theta\cos\varphi, -\sin\theta\sin\varphi, 0) =$$
$$= \frac{1}{\sin\theta}(-\cos\varphi, -\sin\varphi, 0) \qquad 7$$

Hence

$$\vec{u}^2 = \frac{1}{\sin^2\theta}(\sin^2\varphi + \cos^2\varphi) = \frac{1}{\sin^2\theta}$$
$$(\vec{s}\cdot\vec{u}) = \frac{1}{\sin\theta}(-\sin\theta\cos^2\varphi - \sin\theta\sin^2\varphi) = -1 \qquad 8$$

And

$$\Omega = 2\pi - \oint dl \sqrt{\vec{u}^2 - (\vec{s}\cdot\vec{u})^2} = 2\pi - \oint d\varphi \sin\theta \sqrt{\frac{1}{\sin^2\theta} - 1} =$$
$$= 2\pi - \cos\theta \oint d\varphi = 2\pi(1-\cos\theta)$$
(9)

which is a known formula for the solid angle of a cone with apex angle $\theta$.

## 2.3 Derivation

Now let us prove Eq 3 by deriving it from Eq 4. Let us break the smooth parts of the curve into infinitesimally small segments of size $dl$. Let $\delta_j^*$ be an angle enumerated by index $j$ on the sphere between a tangent $\vec{\tau}$ at the point $\vec{s}(l)$ and a vector $\vec{\tau}'$ which is a projection of the tangent $\vec{\tau}(l+dl)$ to the plane normal to $\vec{s}$. Here the smooth curve is approximated by a polygon with $N$ edges and turn angles $\delta_j^*$. As the polygon tends to the curve ($N \to \infty$) all angles $\delta_j^*$ go to zero while $\delta_i$ are fixed at $n$ vertices. The final result can be expressed if the sum of all $\delta_j^*$ and $\delta_i$ between segments is known, Eq 4:

$$\Omega = 2\pi - \sum \delta = 2\pi - \left( \sum_{i=1}^{n} \delta_i + \lim_{N\to\infty} \sum_{j=1}^{N} \delta_j^* \right)$$
(10)

Given that $\vec{\tau}(l+dl) = \vec{\tau} + \vec{u}\,dl$, $(\vec{s}\cdot\vec{\tau}) = 0$, and a projection of any vector $\vec{x}$ on the sphere is $\vec{x} - \vec{s}(\vec{s}\cdot\vec{x})$, $\vec{\tau}'$ becomes:

$$\vec{\tau}' = \vec{\tau}(l+dl) - \vec{s}(\vec{s}\cdot\vec{\tau}(l+dl)) = \vec{\tau} + \vec{u}\,dl - \vec{s}(\vec{s}\cdot(\vec{\tau}+\vec{u}\,dl)) =$$
$$= \vec{\tau} + (\vec{u} - \vec{s}(\vec{s}\cdot\vec{u}))dl$$
(11)

The angles $\delta_j^*$ can be calculated from the scalar products:

$$\cos\delta^* = \frac{\vec{\tau}\cdot\vec{\tau}'}{|\vec{\tau}||\vec{\tau}'|} = \frac{\vec{\tau}\cdot(\vec{\tau} + (\vec{u} - \vec{s}(\vec{s}\cdot\vec{u}))dl)}{|\vec{\tau}|\sqrt{(\vec{\tau}+(\vec{u}-\vec{s}(\vec{s}\cdot\vec{u}))dl)^2}} = \frac{1}{\sqrt{1+dl^2(\vec{u}-\vec{s}(\vec{s}\cdot\vec{u}))^2}}$$
(12)

because $(\vec{u}\cdot\vec{\tau}) = 0$ and $|\vec{\tau}| = 1$. Here index $j$ is implicit. Since both $\delta^*$ and $dl$ are small Eq 12 can be rewritten as:

$$1 - \frac{\delta^{*2}}{2} = 1 - \frac{1}{2}dl^2(\vec{u}-\vec{s}(\vec{s}\cdot\vec{u}))^2 = 1 - \frac{1}{2}dl^2(\vec{u}^2-(\vec{s}\cdot\vec{u})^2) \Rightarrow$$
$$\Rightarrow \delta^* = dl\sqrt{\vec{u}^2 - (\vec{s}\cdot\vec{u})^2}$$
(13)

In the limit the sum of $\delta_j^*$ must be replaced by line integral:

$$\lim_{N\to\infty} \sum_{j=1}^{N} \delta_j^* = \oint dl \sqrt{\vec{u}^2 - (\vec{s}\cdot\vec{u})^2}$$
(14)

Combining together Eqs 10 and 14 proves Eq 3.

## 2.4 Arbitrary smooth generating curve

In a similar line of thought a generic formula for an arbitrary parameterised line in space $L(t)$ without corners can be derived:

$$\Omega = 2\pi - \oint \frac{\left|\vec{L}_1 \times \vec{L}_2\right|}{L_1^2} dt \qquad 15$$

where $L_1$ and $L_2$ are sphere tangential projections of the first and the second derivatives of $L$:

$$\begin{aligned} \vec{L}_1 &= \vec{L}' - \vec{s}(\vec{s} \cdot \vec{L}') \\ \vec{L}_2 &= \vec{L}'' - \vec{s}(\vec{s} \cdot \vec{L}'') \\ \vec{s} &= \vec{L}/L \end{aligned} \qquad 16$$

The numerator in the expression under the integral in Eq 15 is the absolute value of the vector product. Eq 15 reduces to Eq 5 when $L(t) = s(l)$:

$$\begin{aligned} \vec{L}_1 dt &= \vec{\tau}\, dl - \vec{s}(\vec{s} \cdot \vec{\tau})dl = \vec{\tau}\, dl \\ \vec{L}_2 (dt)^2 &= \vec{u}(dl)^2 - \vec{s}(\vec{s} \cdot \vec{u})(dl)^2 \\ (\vec{L}_1 dt)^2 &= (\vec{\tau}\, dl)^2 = (dl)^2 \\ \frac{\left|\vec{L}_1 \times \vec{L}_2\right|}{L_1^2} dt &= \left|\vec{\tau} \times (\vec{u} - \vec{s}(\vec{s} \cdot \vec{u}))\right| dl = \\ &= dl\sqrt{\vec{\tau}^2(\vec{u} - \vec{s}(\vec{s} \cdot \vec{u}))^2 - (\vec{\tau} \cdot (\vec{u} - \vec{s}(\vec{s} \cdot \vec{u})))^2} = dl\sqrt{(\vec{u} - \vec{s}(\vec{s} \cdot \vec{u}))^2} \end{aligned} \qquad 17$$

## 3. Polyhedral cones

### 3.1 Solid angle of discrete shape

Eq 3 is quite general but not always very efficient when many corner angles $\delta$ have to be calculated. Often the *generating curve* is approximated by a sequence of points in which case the line integral of Eq 3 has to be replaced with the sum over corner angles. It might also be important that the sum is calculated in an efficient way.

Let a polyhedral cone be defined by a list of $n$ unit vectors $\{\vec{s}_1, \vec{s}_2, ..., \vec{s}_n\}$. Let $j$ be a circular index in $\vec{s}_j$ so that $\vec{s}_{n+1} = \vec{s}_1$ and $\vec{s}_0 = \vec{s}_n$. Let us also define the following values:

$$\begin{aligned} a_j &= \vec{s}_{j-1} \cdot \vec{s}_{j+1} \\ b_j &= \vec{s}_{j-1} \cdot \vec{s}_j \\ c_j &= \vec{s}_j \cdot \vec{s}_{j+1} = b_{j+1} \\ d_j &= \vec{s}_{j-1} \cdot [\vec{s}_j \times \vec{s}_{j+1}] \end{aligned} \qquad 18$$

Vectors $\{\vec{s}_{j-1}, \vec{s}_j, \vec{s}_{j+1}\}$ form a spherical triangle. *a*, *b*, and *c* are cosines of the triangle sides, and *d* is a volume of the parallelepiped spanned by the vectors $\{\vec{s}_{j-1}, \vec{s}_j, \vec{s}_{j+1}\}$. The angle $\delta_j$ is supplementary to the inward angle $(\pi - \delta_j)$ of the triangle at the vertex $\vec{s}_j$. The *spherical law of cosines* reads (assuming implicit index *j*):

$$a = bc + \sqrt{1-b^2}\sqrt{1-c^2}\cos(\pi-\delta) = bc - \sqrt{1-b^2}\sqrt{1-c^2}\cos\delta \qquad 19$$

from which sine of the angle $\delta$ can be found:

$$d = \sqrt{1-b^2}\sqrt{1-c^2}\sin\delta \qquad 20$$

This can be directly verified by taking into account that $d^2 = 1 + 2abc - a^2 - b^2 - c^2$ (which in turn can be obtained from the determinant formula for a *scalar triple product*).

Knowing sine and cosine of $\delta_j$ the solid angle of the cone can be expressed as

$$\Omega = 2\pi - \sum_{j=1}^{n}\delta_j = 2\pi - \sum_{j=1}^{n}\arctan\frac{d_j}{b_j c_j - a_j} \qquad 21$$

This equation, however, is an inefficient method since it requires taking the arctangent for every corner $\delta_j$. This can be improved by noting that arctangent is the complex argument function of a corresponding complex number $z_j = b_j c_j - a_j + i d_j$:

$$\arctan\frac{d_j}{b_j c_j - a_j} = \arg\{b_j c_j - a_j + i d_j\} = \arg z_j \qquad 22$$

Using the identity:

$$\sum_j \arg z_j = \arg \prod_j z_j \qquad 23$$

which can be proven by

$$\arg\prod_j z_j = \arg\prod_j |z_j| e^{i\arg z_j} = \arg\prod_j e^{i\arg z_j} = \arg e^{i\sum_j \arg z_j} = \sum_j \arg z_j \qquad 24$$

Eq 21 becomes:

$$\boxed{\Omega = 2\pi - \arg\prod_{j=1}^{n}\{b_j c_j - a_j + i d_j\}} \qquad 25$$

This formula is of *O(n)* complexity requiring 3*n* scalar (note that $b_{j+1} = c_j$) and *n* vector multiplications of vectors, *n* multiplications of complex numbers, and only one arctangent. In contrast, a derived formula in Ref [2] requires square root and trigonometric functions for each vertex.

## 3.2 Tetrahedron example

A nice formula is derived for the solid angle of a tetrahedron in Ref [3]. It can be reproduced using the result of Eq 25. In a triangular case $d_j = d$ are the same and the product in Eq 25 is expressed via sides $\{a,b,c\}$ of the spherical triangle and can be factored in the following way:

$$\Omega = 2\pi - \arg(bc - a + id)(ac - b + id)(ab - c + id) =$$
$$= 2\pi - \arg \frac{1}{2}(1-a)(1-b)(1-c)(1+a+b+c-id)^2 = \qquad 26$$
$$= 2\pi - \arg(1+a+b+c-id)^2 = 2\pi - 2\arg(1+a+b+c-id)$$

Here the expression for $d^2$ has been used as well as the following two identities for $arg$ function: $\arg z^n = n \arg z$ and $\arg kz = \arg z$ for any $k > 0$.

Dividing Eq 26 by 2 and taking tangent one finds:

$$\tan \frac{\Omega}{2} = \tan(\pi - \arg(1+a+b+c-id)) = \tan\arg(1+a+b+c+id) =$$
$$= \frac{d}{1+a+b+c} \qquad 27$$

which is the formula for spherical triangle derived in Ref [3].

The symmetry of the products in Eq 26 allows factoring the expressions in the triangular case and reducing it to a simple ratio as Eq 27. It is difficult to say whether it is possible to simplify Eq 25 in a similar manner for the general $n$-polygonal case.

## 3.3 Face of a cube example

As another example consider the solid angle subtended by one face of a cube of side length 2 centred at the origin. Since the cube is symmetrical and has six sides, one side obviously subtends angle $4\pi/6=2\pi/3$. Vertices of the cube can be placed at points: $(1, 1, 1)$, $(1, -1, 1)$, $(-1, 1, 1)$, and so on. Since the corners are the same in the cube, values $a$, $b$, $c$, and $d$ are the same for each $j$. Hence only one spherical triangle has to be calculated. Normalising vectors to unit:

$$\bar{s}_1 = \frac{1}{\sqrt{3}}(1, -1, 1)$$
$$\bar{s}_2 = \frac{1}{\sqrt{3}}(1, 1, 1) \qquad 28$$
$$\bar{s}_3 = \frac{1}{\sqrt{3}}(-1, 1, 1)$$

gives

$$a = (\vec{s}_1 \cdot \vec{s}_3) = -\frac{1}{3}$$

$$b = (\vec{s}_1 \cdot \vec{s}_2) = \frac{1}{3}$$

$$c_j = (\vec{s}_2 \cdot \vec{s}_3) = \frac{1}{3} \tag{29}$$

$$d_j = (\vec{s}_1 \cdot [\vec{s}_2 \times \vec{s}_3]) = \frac{4}{3\sqrt{3}}$$

Then

$$\Omega = 2\pi - \arg \prod_{j=1}^{4} \{bc - a + id\} = 2\pi - \arg\left(\frac{1}{9} + \frac{1}{3} + i\frac{4}{3\sqrt{3}}\right)^4 =$$
$$= 2\pi - 4\arg(1 + i\sqrt{3}) = 2\pi - 4\arctan\sqrt{3} = 2\pi - 4\pi/3 = 2\pi/3 \tag{30}$$

as expected.

## 4. Intersecting cones

Imagine two cones with the same apex. Each cone has its own axis and apex angle. Depending on apex angles and the angle between the cone's axes, one cone can be entirely inside the other, entirely outside, or partially inside and outside. In the later case surfaces of the cones are intersecting. It is easy to answer a question about solid angle of both cones when one is completely inside or outside the other. When the cones are intersecting the total solid angle of both cones is equal to the sum of each cone independently minus the solid angle of their intersection. Hence the difficult question is to find their intersection.

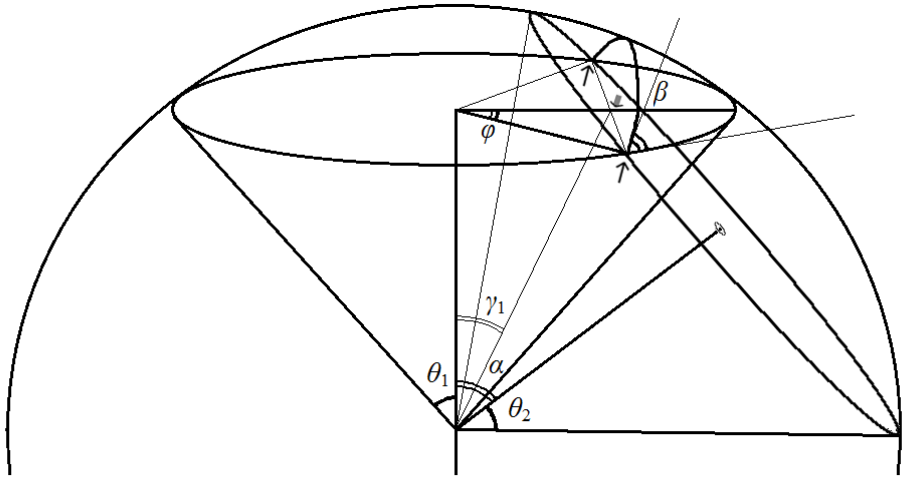

**Figure 1** Two intersecting cones inside a sphere. *Intersection points* are shown with two small arrows. The thick short arrow shows *Middle line*, which passes through the origin and the centre of the line segment between *Intersection points*. *Intersection points* are connected with three arcs: two arcs of *Circle*s and one arc of Great circle. The angle *β* is an angle between tangents of *Great circle* arc and the arc of one of the *Circle*s.

Considering Figure 1 let us define the following entities.

*Cone* is one of two cones sharing the same apex
*Sphere* is an imaginary unit sphere centred at apex of the cones
*Circle* is a circle made by intersection of *Cone*'s surface and *Sphere*
*Plane* is a plane constructed by two lines of intersection of *Cone*s surfaces
*Great circle* is a circle made by intersection of *Plane* and *Sphere*
*Intersection points* are points where two *Circles* of two *Cones* intersect
*Cap* is a part of *Sphere*'s surface cut by *Cone*
*Cap's segment* is a part of *Cap* cut by *Plane*
*Middle line* is a line made by the intersection of *Plane* and the plane of *Cone*s axes.

*β* – angle of intersection between *Circle* and *Great circle*
*φ* – half-angle at *Circle*'s centre between *Intersection points*
*γ* – angle between *Plane* and *Cone*'s axis
*α* – angle between two *Cone*s axes
*θ* – apex angle of *Cone*

It is obvious that the solid angle of intersection of two cones is equal to the sum of the solid angles of two *Cone*'s segments[1]. Hence let us first find the solid angle of *Cap*'s *segment* – the part of a *Cone* cut by a *Plane* at angle *γ*.

## 4.1 Cone's Segment

The solid angle of *Cone*'s segment is determined by two angles *θ* and *γ*: $\Omega = \Omega(\theta, \gamma)$.

### 4.1.1 Closed curve on Sphere

The curve defining *Cap's segment* consists of Circle's arc of 2*φ*, two turn angles $(\pi - \beta)$, and the arc of *Great circle*. Using Eq 3 and noting that the arc of *Great circle* does not contribute to the formula one finds:

$$\Omega = 2\pi - (\pi - \beta) - (\pi - \beta) - \int_{-\varphi}^{\varphi} d\varphi' \cos\theta = 2(\beta - \varphi \cos\theta) \qquad 31$$

Now it is a matter to calculate angles *β* and *φ*.

### 4.1.2 Circle sector angle

Let *Cone*'s axis be aligned with Z-axis and *Circle* be defined as:

$$Circle = (\sin\theta\cos\varphi, \sin\theta\sin\varphi, \cos\theta) \qquad 32$$

Let *Great circle* be defined as

$$GreatCircle = (\sin\gamma\cos\varphi', \sin\varphi', \cos\gamma\cos\varphi') \qquad 33$$

which is easily obtained rotating a vertical great circle $(0, \sin\varphi', \cos\varphi')$ around the Y-axis by angle *γ*. Here *γ* is the angle between *Plane* and Z-axis. $\varphi'$ and *φ* are

---

[1] *Cone*'s segment is a solid and *Cap's segment* is a surface, but the solid angle of each one refers to the same entity.

~ 8 ~

different angles each defined within its own circle. Comparing X and Z components of circle *Intersection points* gives the result:

$$\begin{cases} \sin\theta\cos\varphi = \sin\gamma\cos\varphi' \\ \cos\theta = \cos\gamma\cos\varphi' \end{cases} \Rightarrow \tan\theta\cos\varphi = \tan\gamma \qquad 34$$

$\varphi$ is then defined by the expression

$$\cos\varphi = \frac{\tan\gamma}{\tan\theta} \qquad 35$$

### 4.1.3 Tangent angle

Angle $\beta$ can be found by the formula:

$$\cos\beta = \frac{\vec{\tau}_1 \cdot \vec{\tau}_2}{\tau_1 \tau_2} \qquad 36$$

where $\tau$ are tangent vectors at *Intersection points* which in turn are found by taking derivative of Eq 32 and Eq 33:

$$\begin{aligned} \vec{\tau}_1 &= (-\sin\theta\sin\varphi, \sin\theta\cos\varphi, 0) \\ \vec{\tau}_2 &= (-\sin\gamma\sin\varphi', \cos\varphi', -\cos\gamma\sin\varphi') = \\ &= \left(-\sin\gamma\sin\theta\sin\varphi, \frac{\cos\theta}{\cos\gamma}, -\cos\gamma\sin\theta\sin\varphi\right) \end{aligned} \qquad 37$$

Next is to find the ingredients for Eq 36:

$$\begin{aligned} \tau_1 &= \sqrt{(\sin\theta\sin\varphi)^2 + (\sin\theta\cos\varphi)^2} = \sin\theta \\ \tau_2 &= 1 \\ \vec{\tau}_1 \cdot \vec{\tau}_2 &= \sin\gamma\sin^2\theta\sin^2\varphi + \sin\theta\cos\varphi\frac{\cos\theta}{\cos\gamma} = \sin\gamma \end{aligned} \qquad 38$$

The last step in the above equation is a trigonometric exercise taking into account Eq 35. Inserting Eq 38 into Eq 36 produces[2]:

$$\cos\beta = \frac{\sin\gamma}{\sin\theta} \qquad 39$$

Note that Eqs 35 and 39 bind the angles by the expression:

---

[2] When deriving Eq 39 I obtained a result: $\tan\beta\tan\theta = \tan\varphi$. Later I found a mistake in the derivation. I could not believe that the result would be an uglier formula. Fortunately the correct result is as nice which made me wonder what the corresponding interpretation of the wrong result was. It turned out this result corresponds to the case of cylinder instead of sphere. The mistake was that the great circle was defined as $(\sin\gamma\cos\varphi', \sin\gamma\sin\varphi', \cos\gamma\cos\varphi')$.

$$\cos\varphi\cos\gamma = \cos\beta\cos\theta \qquad 40$$

### 4.1.4 Solid angle of the segment

Combining Eqs 31, 35, and 39 one gets the final formula for the solid angle of *Cone*'s segment:

$$\boxed{\Omega(\theta,\gamma) = 2\left(\arccos\frac{\sin\gamma}{\sin\theta} - \cos\theta\cdot\arccos\frac{\tan\gamma}{\tan\theta}\right)} \qquad 41$$

## 4.2 Intersecting plane

Once we know how to calculate the solid angle of Cone's segment, finding the intersection of two cones requires finding the *Plane*'s angle $\gamma$.

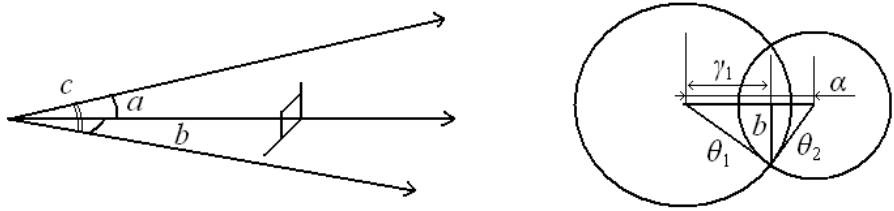

**Figure 2** The left picture shows a right spherical triangle. The right picture show *Circles* projected to a Euclidian plane and spherical angles between *Circle* centres, *Middle line* and one of *Intersection points*.

First note the centres of *Caps*, *Middle line*, and one of the *Intersection points* form two right spherical triangles with sides $(\gamma_1,\theta_1,b)$ and $(\gamma_2,\theta_2,b)$ in Figure 2 right. $b$ is the common cathetus and $\gamma_1 + \gamma_2 = \alpha$. A right spherical triangle $(a,b,c)$, see Figure 2 left, where $c$ is hypotenuse can be constructed by three vectors $(1,0,0)$, $(\cos a, 0, \sin a)$, and $(\cos b, \sin b, 0)$. Cosine of $c$ is a scalar product of two later vectors:

$$\cos c = (\cos a, 0, \sin a)\cdot(\cos b, \sin b, 0) = \cos a \cos b \qquad 42$$

From here

$$\begin{cases} \cos\theta_1 = \cos\gamma_1 \cos b \\ \cos\theta_2 = \cos\gamma_2 \cos b \end{cases} \Rightarrow \frac{\cos\theta_2}{\cos\theta_1} = \frac{\cos(\alpha-\gamma_1)}{\cos\gamma_1} = \cos\alpha + \sin\alpha\tan\gamma_1 \qquad 43$$

which gives:

$$\qquad 44$$

$$\boxed{\tan \gamma_1 = \frac{\cos \theta_2 - \cos \alpha \cos \theta_1}{\sin \alpha \cos \theta_1}}$$

## 4.3 The solution

The solid angle of intersection of two cones defined by three numbers $\theta_1$, $\theta_2$ and $\alpha$ can be calculated by the following set of equations:

$$\boxed{\begin{aligned} \Omega &= \Omega(\theta_1, \gamma_1) + \Omega(\theta_2, \gamma_2) \\ \Omega(\theta, \gamma) &= 2(\beta - \varphi \cos \theta) \\ \cos \varphi &= \frac{\tan \gamma}{\tan \theta} \\ \cos \beta &= \frac{\sin \gamma}{\sin \theta} \\ \tan \gamma_1 &= \frac{\cos \theta_2 - \cos \alpha \cos \theta_1}{\sin \alpha \cos \theta_1} \end{aligned}} \qquad 45$$

$\gamma_2$ is calculated by exchanging indices 1 and 2. Ranges of the parameters are:

$$\begin{aligned} 0 &< \alpha < \pi \\ 0 &< \theta_{1,2} < \pi/2 \end{aligned} \qquad 46$$

For the cases $\pi/2 < \theta_{1,2} < \pi$ see Section 4.4.2 below.

## 4.4 Numerical considerations

This section describes practical recommendations for avoiding singularities in the solution formulae. They are presented in the order of generality, so formulae do not break as long as the previous cases were taken care of.

The solid angle of each cone used below in this Section is expressed via its apex angle by a known formula (derived also as an example in Eqs 6-9):

$$\Omega_{1,2} = 2\pi(1 - \cos \theta_{1,2}) \qquad 47$$

If a subsection below mentions only $Cone_1$, the same argument applies to $Cone_2$ and vice versa.

### 4.4.1 Avoiding Arctan and illegal Cos

To avoid Arctan in calculation $\gamma$ two values can be calculated:

$$\begin{aligned} t_y &= \cos \theta_2 - \cos \alpha \cos \theta_1 \\ t_x &= \sin \alpha \cos \theta_1 \end{aligned} \qquad 48$$

Now $\varphi$ and $\beta$ can be calculated as:

$$\cos\varphi_1 = \left[\frac{t_y \cos\theta_1}{t_x \sin\theta_1}\right]_{-1}^{+1} \qquad \textbf{49}$$

$$\cos\beta_1 = \left[\frac{t_y}{\sin\theta_1 \sqrt{t_x^2 + t_y^2}}\right]_{-1}^{+1} \qquad \textbf{50}$$

where [] is a limiting function:

$$[x]_{-1}^{+1} = \begin{cases} 1 & \text{if} \quad x > 1 \\ -1 & \text{if} \quad x < -1 \\ x & \text{otherwise} \end{cases} \qquad \textbf{51}$$

The limiting function appears because the cases when cosines of $\varphi$ and $\beta$ are outside the range (-1,1) correspond to the situations when one cone is completely inside or outside the other.

### 4.4.2 Inverted cones

If the apex angle of $Cone_1$ is greater than $\pi/2$, then the correct result can be obtained by the following substitution:

$$\Omega(\theta_1 > \pi/2) = \Omega_2 - \Omega(\alpha \to \pi - \alpha, \theta_1 \to \pi - \theta_1) \qquad \textbf{52}$$

### 4.4.3 Co-directed cones

When $\alpha \to 0$ then the result approaches to the smaller cone being inside the other:

$$\Omega = \min(\Omega_1, \Omega_2) \qquad \textbf{53}$$

### 4.4.4 Counter-directed cones

When $\alpha \to \pi$ then the result approaches to an overlap between $Cone$s or zero if their sum is less than the solid angle of a sphere:

$$\Omega = \begin{cases} \Omega_1 + \Omega_2 - 4\pi & \text{if} \quad > 0 \\ 0 & \text{otherwise} \end{cases} \qquad \textbf{54}$$

### 4.4.5 Narrow cone

When $Cone_1$ is narrow, $\theta_1 \to 0$, then $\gamma_1 = \alpha - \theta_2$. Excluding the case when $\theta_1 = 0$ exactly, the result for simplicity can be obtained by linear approximation:

$$\Omega = \begin{cases} 0 & \text{if} \quad \gamma_1 > \theta_1 \\ \Omega_1 & \text{if} \quad \gamma_1 < -\theta_1 \\ \Omega_1 \dfrac{\gamma_1 + \theta_1}{2\theta_1} & \text{otherwise} \end{cases} \qquad \textbf{55}$$

### 4.4.6 Two hemispheres

Another special case is when both *Cone*s approach to hemispheres $\theta_1 \to \pi/2$ and $\theta_2 \to \pi/2$. This case is different (and has to be checked first) from the case when just one cone becomes a hemisphere. The intersection is linear to $\alpha$ changing from $2\pi$ when $\alpha = 0$ to 0 when $\alpha = \pi$:

$$\Omega = 2(\pi - \alpha) \qquad 56$$

### 4.4.7 One hemisphere

In the case of one hemisphere, $\theta_1 \to \pi/2$ the intersection is equal to the segment of the spherical cap of the *Cone₂*:

$$\Omega = \Omega(\theta_2, \gamma_2) = \Omega(\theta_2, \alpha - \pi/2) \qquad 57$$

## 4.5 Monte-Carlo simulation

Figure 3 shows an example of intersecting *Cone*s with $\Omega_1 = 0.6 \cdot 4\pi$ and $\Omega_2 = 0.2 \cdot 4\pi$ checked with Monte-Carlo simulations calculating solid angle by putting random points on a unit sphere surface. The right side graph shows how the exact solution deviates from the linear approximation.

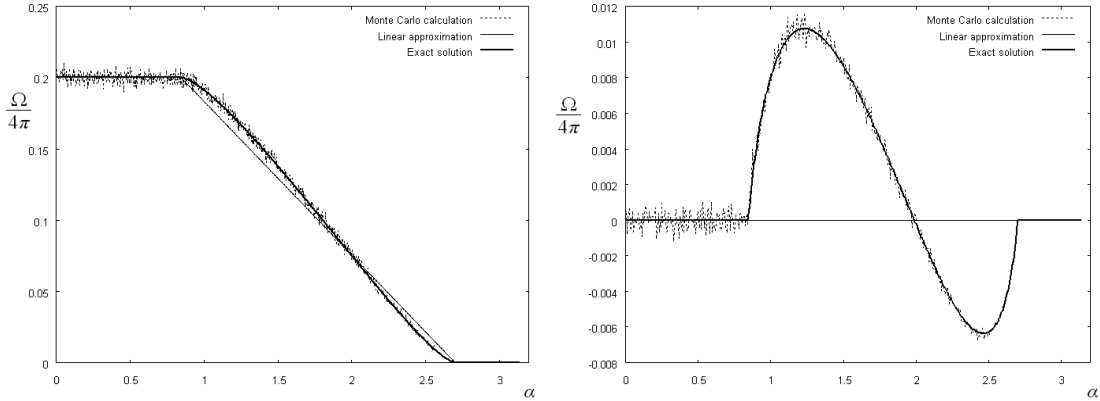

**Figure 3** Example of two intersecting cones. The X-axis is the angle between axes of the *Cone*s. The Y-axis is the intersection solid angle normalised to $4\pi$. The thin solid line is the exact solution when one *Cone* is entirely inside or outside the other and linear approximation between those two cases. The thick line is the solution, Eq 45. The dashed line is Monte-Carlo simulation. The right picture shows the difference between a corresponding curve and the linear approximation. The left and the right pictures have different Monte-Carlo statistic.

# 5. Conclusion

In this paper I investigated some aspects of calculating the solid angles of conical shapes. A solution for the solid angle of two intersecting cones is obtained in Eq 45. Brief derivation of this result and advice for numerical computation are covered. The solid angles of arbitrary conical surfaces defined by a closed curve on a unit sphere, Eq 3, and a more general expression for arbitrary smooth closed curves, Eq 15, are presented. A new efficient formula for polyhedral cones, Eq 25, is derived.